\input xy
\xyoption{all}
\input amssym.tex

\mag=1200 
\hsize=130mm  \vsize=180mm  \voffset 5mm
\lineskiplimit=10pt \lineskip=10pt
\tolerance=10000 \pretolerance=1000 \parindent=0mm \raggedright

 \font\ud =cmr10 at  12pt 
 \font\uq =cmr10 at  14pt

 \font\dz =cmr10 at  20pt

\font\tengo=eufm10

\font\sevengo=eufm7

\font\fivego=eufm5

\font\tenbb=msbm7 at 10pt

\font\sevenbb=msbm7   
 
\font\fivebb=msbm5


\newfam\gofam  \textfont\gofam=\tengo
\scriptfont\gofam=\sevengo   \scriptscriptfont\gofam=\fivego
\def\go{\fam\gofam\tengo}

\newfam\bbfam  \textfont\bbfam=\tenbb
\scriptfont\bbfam=\sevenbb   \scriptscriptfont\bbfam=\fivebb
\def\bb{\fam\bbfam\tenbb}


\def\ind{\hskip 1em\relax}

\def\Hom{\mathop{\rm Hom}\nolimits}

\let\f=\varphi
\let\to=\rightarrow

\let\mbox=\hbox

\def\Ext{\mathop{\rm Ext}\nolimits}
 \def\C{{\bb C}}  \def\Z{{\bb Z}}

\def\cl{\centerline}

\def\{{\lbrace}  
\def\}{\rbrace}   
\def\({\langle}  
\def\){\rangle}
\def\[{\lbrack} \def\]{\rbrack}

\def\arrow{\rightarrow}
\def\iso{\buildrel\sim\over{\arrow}}
\def\.{\bullet}
\def\bs{\bigskip}

\def\mono{\rightarrowtail}
\def\epi{\twoheadrightarrow}
\def\incl{\hookrightarrow}

\def\wt{\widetilde}

\def\t{\mathop{\otimes}}  
\def\ds{\mathop{\oplus}}

\def\\{\backslash}
\def\back/{\backslash}
\def\e{\varepsilon}
%
\def\hfl#1#2{\smash{\mathop{\hbox to 12mm{\rightarrowfill}}
\limits^{\scriptstyle#1}_{\scriptstyle#2}}}
\def\vfl#1#2{\llap{$\scriptstyle#1$}\left\downarrow
\vbox to 6mm{}\right.\rlap{$\scriptstyle#2$}}
\def\hfle#1#2{\smash{\mathop{\hbox to 12mm{\leftarrowfill}}
\limits^{\scriptstyle#1}_{\scriptstyle#2}}}
\def\vfle#1#2{\llap{$\scriptstyle#1$}\left\uparrow
\vbox to 6mm{}\right.\rlap{$\scriptstyle#2$}}
\def\hfld#1#2{\smash{\mathop
{\buildrel {\hbox to 12mm{\rightarrowfill}} 
\over{{\hbox to 12mm{\leftarrowfill}}}}
\limits^{\scriptstyle#1}_{\scriptstyle#2}}} 
\def\diagram#1{\def\normalbaselines{\baselineskip=0pt
\lineskip=10pt\lineskiplimit=1pt} \matrix{#1}}
%
%
\cl{\dz Four \ Multiplicative \ Cohomology Theorems} \bs
The main statement proved here is the \bs
\cl{\uq De \ \ Rham \ \ Theorem} \bs
{\ud The de Rham cohomology algebra of a paracompact manifold
is canonically isomorphic to the graded algebra of endomorphisms
of the constant sheaf viewed as an object in the derived category
of sheaves of vector spaces. } \bs
There are similar statements for Lie algebra cohomology and group cohomology.  
\bs \hrule \bs
\cl{\ud Reminder \ about \ the \ notion \ of \ derived \ category} \bs
Let $\cal A$ be an abelian category and ${\cal A}^\.$ the category of cochain complexes 
in $\cal A$. Recall that a morphism in ${\cal A}^\.$ is a {\bf quasi-isomorphism} if it 
induces an isomorphism in cohomology. 
A {\bf derived category} of $\cal A$ is a pair $({\cal D},i)$ where $\cal D$ 
is a category and $i$ a functor from ${\cal A}^\.$ to $\cal D$ satisfying the 
following condition. \bs
\ind For any category $\cal C$ and any functor $F : {\cal A}^\. \arrow {\cal C}$ 
transforming  quasi-isomorphisms into isomorphisms there is a unique functor 
$\wt{F} : {\cal D \arrow \cal C}$ such that

$$\xymatrix{
{\cal A}^\. \ar[r]^F \ar[d]_i & {\cal C} \\
{\cal D} \ar@{-->}[ur]_{\wt{F}}}$$
commutes. \bs
\ind It's easy to check that a derived category of $\cal A$ exists and is 
unique 
up to unique isomorphism, and that $\cal D$ has the same objects as 
${\cal A}^\.$. 
It can be shown that $\cal D$ is additive. \bs
\ind For $C^\. \in {\cal A}^\.$ and $n \in \Z$ denote by $C^{n+\.}$ the 
indicated shifted 
complex equipped with $(-1)^n$ times the differential of $C^\.$. For 
$C^\., D^\. \in {\cal A}^\.$ and $n \in \Z$ put
$$\Ext_{\cal A}^n(C^\., D^\.) := \Hom_{\cal D}(C^\., D^{n+\.})$$
and note that there is an obvious composition
$$\Ext_{\cal A}^q(D^\., E^\.) \t_\Z \Ext_{\cal A}^p(C^\., D^\.) 
\longrightarrow
\Ext_{\cal A}^{p+q}(C^\., E^\.),$$
called {\bf cup-product}.
Sometimes one abusively calls {\bf derived category} the category $\cal E$ 
whose objects are those of $\cal A^\.$ and whose morphisms are defined by 
the rule

$$\Hom_{\cal E}(C^\.,D^\.) := 
\bigoplus_{n \in \Z} \ \Ext_{\cal A}^n(C^\., D^\.).$$ \hrule \bs
This short text is about the four following cohomology theories : the de 
Rham cohomology, the (relative) Lie algebra cohomology, the group cohomology 
and the \v Cech cohmology (I refer to Warner [W], Borel-Wallach [BW] and 
Cartan-Eilenberg [CE] for precise definitions, and to Verdier [V] for derived 
category theory). Each of these theories was first defined by a magic 
formula and then interpreted as an Ext-group~;  in each case there is an 
obvious cup-product suggested by the formula and an obvious cup-product 
suggested by the Ext interpretation~; the goal is to show that the 
combinatorial cup-product is compatible with the conceptual one. \eject
\cl{\uq Statement \ of \ the \ Theorems} \bs
Let $M$ be a paracompact  manifold, $\C_M$ the constant 
sheaf with fiber $\C$ over $M$ and $\C_M \mono \Omega$ the de Rham 
resolution. \bs
{\bf Multiplicative de Rham Theorem.} The de Rham 
cohomology of $M$ is canonically isomorphic as a graded algebra to 
$\Ext^\._{\C_M\hbox{\sevenrm -mod}}(\C_M,\C_M)$. \bs
Let $k$ be a field of characteristic 0, let \hbox{$\go k \subset \go g$}  
be finite dimensional Lie algebras over $k$~; assume $\go g$ is $\go k$-semisimple. 
Given a $(\go g,k)$-module ({\it i.e.} a $\go k$-semisimple \hbox{$\go g$-module)} 
$V$, denote the 
Chevalley-Eilenberg cohomology of $({\go g,k})$ with values in $V$ by 
$H_{\hbox{\sevenrm CE}}^\.({\go g,k} \, ; V)$ and put 
$H^\.({\go g,k} \, ; V) := \Ext_{\go g,k}^\.(k,V).$ By a Theorem of Hochschild there is
a canonical isomorphism of graded vector spaces 
$$\Phi_{V}: H_{\hbox{\sevenrm CE}}^\.({\go g,k} \, ; V) \iso H^\.({\go g,k} \, ; V).$$
Consider the diagram 
$$\diagram{
H_{\hbox{\sevenrm CE}}^\.({\go g,k} \, ; V) \t 
H_{\hbox{\sevenrm CE}}^\.({\go g,k} \, ; k) & 
\hfl{ \ \Phi_V \t \Phi_k \ \ }{\sim} & 
H^\.({\go g,k} \, ; V) \t H^\.({\go g,k} \, ; k) \cr
\vfl{}{} &                    & \vfl{}{} \cr
H_{\hbox{\sevenrm CE}}^\.({\go g,k} \, ; V)  & \hfl{\sim}{\Phi_V} & 
H^\.({\go g,k} \, ; V), \cr}$$
where the vertical arrows represent the cup-products. \bs 
{\bf Multiplicative Hochschild Theorem.} This diagram commutes. \bs
Let $G$ be a group and $k$ a commutative ring. Denote the 
Eilenberg-MacLane cohomology of $G$ with values in the $G$-module $V$ by 
$H_{\hbox{\sevenrm EM}}^\.(G,V)$ and put 
$H^\.(G,V) := \Ext_G^\.(k,V).$ By a Theorem of Eilenberg-MacLane there is
a canonical isomorphism of graded $k$-modules
$$\Phi_{V}: H_{\hbox{\sevenrm EM}}^\.(G,V) \iso H^\.(G,V).$$
Consider the diagram 
$$\diagram{
H_{\hbox{\sevenrm EM}}^\.(G,V) \t 
H_{\hbox{\sevenrm EM}}^\.(G,k) & 
\hfl{ \ \Phi_V \t \Phi_k \ \ }{\sim} & 
H^\.(G,V) \t H^\.(G,k) \cr
\vfl{}{} &                    & \vfl{}{} \cr
H_{\hbox{\sevenrm EM}}^\.(G,V)  & \hfl{\sim}{\Phi_V} & 
H^\.(G,V), \cr}$$
where the vertical arrows represent the cup-products. \bs 
\vbox{{\bf Multiplicative Eilenberg-MacLane Theorem.} This diagram commutes.} \bs
Given an  be an open cover $\cal U$ of a topological space $X$, a commutative ring $k$ and a sheaf $S$ of modules over the constant sheaf $C := k_X$ there is --- by a results I'll call \v Cech Theorem --- a canonical morphism of graded $k$-modules
$$\Phi_{S}: \check{H}^\.({\cal U},S) \to H^\.(X,S).$$
Consider the diagram 
$$\diagram{
\check{H}^\.({\cal U},S) \t 
\check{H}^\.({\cal U},C) & 
\hfl{ \ \Phi_S \t \Phi_C \ \ }{} & 
H^\.(X,S) \t H^\.(X,C) \cr
\vfl{c}{} &                    & \vfl{}{t} \cr
\check{H}^\.({\cal U},S)  & \hfl{}{\Phi_S} & 
H^\.(X,S), \cr}$$
where $c$ denotes the \v Cech cup-product and $t$ denotes the ``true'' 
cup-product. \bs 
{\bf Multiplicative \v Cech Theorem.} This diagram commutes. \bs
Here is a theorem that is not multiplicative at all~; I put it here because 
its proof uses \v Cech cohomology --- and because I don't know where else to 
put it. \bs
\vbox{{\bf Bott-Tu Theorem.} The Fr\'echet structure of the de Rham cohomology 
of a paracompact manifold is a topological invariant.} \bs
{\bf Proof.} In Proposition 9.8 of [BT] Bott and Tu give an explicit 
isomorphism $\f$ from the cohomology $\check{H}^\.({\cal U},\C)$ of a countable good 
cover onto the de Rham cohomology $H^\._{\sevenrm dR}(M,\C)$. A look at their 
formula shows that $\f$ is continuous. Since the target is Fr\'echet, the 
Open Mapping Theorem implies that $\f$ is a Fr\'echet isomorphism. Therefore 
the Fr\'echet structure of $\check{H}^\.({\cal U},\C)$ doesn't depend on the choice of 
$\cal U$, and, similarly, the Fr\'echet structure of 
$H^\._{\sevenrm dR}(M,\C)$ 
doesn't depend on the choice of the differential structure of $M$. This 
proves that the Fr\'echet structure in question depends only on the topology 
of $M$. {\bf QED} \bs
%
%
\cl{\uq Statement \ of \ the \ Propositions} \bs
Let $k$ be a commutative ring and $\cal C$ a $k$-category~; the convention 
$$\t := \t_k$$
shall be in force to the end. The following abbreviations will come 
handy : if $A$ and $B$ are complexes in $\cal C$ then $\(A,B \)$ shall be the 
complex whose underlying graded $k$-module is defined by 
$$\(A,B \)^n \ := \ \ds_p \ \Hom(C^p,D^{n+p}),$$
the differential being given by $df = d \circ f - (-1)^n \ f \circ d$ for 
$f \in \(A,B\)^n$, and $[A,B]$ shall be its cohomology : 
$$[A,B] := H(\(A,B \)),$$
recalling that $[A,B]^n$ is also the \hbox{$k$-module} of 
homotopy classes of complex morphisms from $A^\.$ to $B^{n+\.}$. 
Let 
$$can_{AB} : [A,B] \arrow \Ext(A,B)$$
be the natural morphism. If $C$ is a third complex then we have the 
composition morphisms 
$$\(B,C \) \ \t \ \(A,B \) \arrow \(A,C \),$$
$$[B,C] \ \t \ [A,B] \arrow [A,C],$$
$$\Ext(B,C) \ \t \ \Ext(A,B) \arrow \Ext(A,C),$$
all abusively denoted by $c$ and the last one being by definition the cup-product~; 
moreover the following diagram commutes
$$\diagram{
[B,C] \ \t \ [A,B] & \hfl{c}{} & [A,C]  \cr 
\vfl{can_{BC} \t can_{AB}}{} &           & \vfl{}{can_{AC}} \cr 
\Ext(B,C) \ \t \ \Ext(A,B) & \hfl{}{c}  & \Ext(A,C). \cr} \eqno{(1)} $$
\ind For $i=1,2,3$ let $V_i$ be an object 
of $\cal C$ and $\e_i : V_i \mono A_i$ a right resolution of $V_i$. Recall that 
we have the commuting diagram of isomorphisms
$$\diagram{
\Ext(A_i,V_j) & \hfl{\Ext(\e_i,V_j)}{} & \Ext(V_i,V_j)  \cr 
\vfl{\Ext(A_i, \e_j)}{} &           & \vfl{}{\Ext(V_i,\e_j)} \cr 
\Ext(A_i,A_j) & \hfl{}{\Ext(\e_i,A_j)}  & \Ext(V_i,A_j), \cr} $$ 
which I'll implicitly use to identify these four $k$-modules whenever convenient. 
Also I'll often write $f \, g$ for $f \circ g$.  \bs
\ind Assume we are given complex morphisms 
$$\f_{ij} : \( V_i,A_j \) \arrow \( A_i,A_j \)$$
for $i < j$, and 
$$\mu : \( V_2,A_3 \) \ \t \ \(V_1,A_2 \) \arrow \( V_1,A_3 \)$$
subject to \bs
{\parindent=7mm
\item{(a)}  $\f_{ij}(f) \ \e_i = f$ for all $f$ in $\(V_i,A_j\)$ and \bs
\item{(b)} the diagram
$$\diagram{
\( V_2,A_3 \) \ \t \ \(V_1,A_2 \) & \hfl{\mu}{} & \( V_1,A_3 \)  \cr 
\vfl{ \f_{{}_{23}} \t \f_{{}_{12}}}{} & & \vfl{}{\f_{{}_{13}}} \cr 
\( A_2,A_3 \) \ \t \ \(A_1,A_2 \) & \hfl{}{c}  & \( A_1,A_3 \) \cr} $$ 
\item{} commutes up to homotopy. \par} \bs
Let $\mu_* : [V_2,A_3] \t \[V_1,A_2] \arrow [V_1,A_3]$ be the 
morphism induced by $\mu$. Denote both $can_{{}_{V_i A_j}}$ and $can_{{}_{A_i A_j}}$
by $can_{ij}$ when convenient. --- Consider the diagram
$$\diagram{
[V_2,A_3] \ \t \ [V_1,A_2] \ & \hfl{\mu_*}{} & [V_1,A_3]  \cr 
\vfl{can_{23} \t can_{12}}{}           &             & \vfl{}{can_{13}} \cr 
\Ext(V_2,V_3) \ \t \ \Ext(V_1,V_2)     & \hfl{}{c}   & \Ext(V_1,V_3). \cr} $$ 
{\bf Proposition 1}. This diagram commutes. \bs
{\bf Weak Proposition 2}. If $\cal C$ has enough injectives and 
\hbox{$\Ext^p(V_i,A_j^q) = 0$} for \hbox{$p > 0$}, \hbox{$q \geq 0$} then 
$can_{ij} : [V_i,A_j] \arrow \Ext(V_i,V_j)$ is an isomorphism. \bs
{\bf Proof}. This follows immediately from Grothendieck's Remark 3 after 
Theorem 2.4.1 in [G]. {\bf QED} \bs
\vbox{{\bf Strong Proposition 2}.  If \hbox{$\Ext^p(V_i,A_j^q) = 0$} for \hbox{$p > 0$}, 
\hbox{$q \geq 0$} then 
$$can_{ij} : [V_i,A_j] \arrow \Ext(V_i,V_j)$$
is an isomorphism.} \bs 
\vbox{\cl{\uq Proposition \ 1 \ and \ Weak \ Proposition \ 2} \bs
\cl{\uq imply \ the \ Theorems} \bs
{\bf Proof of the Multiplicative de Rham Theorem}. Let 
$\C_M \mono \Omega$ the de Rham resolution. 
Put \hbox{$\Omega(M) := \( \C_M,\Omega \)$} as usual, and 
\hbox{$V_1 = V_2 = V_3 := \C_M$,} $A_1 = A_2 = A_3 := \Omega$, and define 
$\f : \Omega(M) \arrow \(\Omega,\Omega \)$ by 
$\f(\alpha)(\xi) = \alpha \wedge \xi$ (here $\alpha$ is a form and $\xi$ a germ) and 
$\mu : \Omega(M) \t \Omega(M) \arrow \Omega(M)$ by 
$\mu(\alpha,\beta) := \alpha \wedge \beta$. {\bf QED}} \bs
{\bf Proof of the Multiplicative Hochschild Theorem}. 
Let $\cal C$ be the category of \hbox{$(\go g,k)$-modules.} 
Note that if $W$ is a complex in $\cal C$ then $\( k , W \)$ is the subcomplex  
$W^{\go g}$ of invariants. Put \hbox{$V_1 = V_2 := k$}, $V_3 := V$, 
$$A_i := \Hom_k \Big( \ U({\go g}) \ \t_{\go k} \ \bigwedge({\go g/k}) \ , \ V_i \ 
\Big)_{K \hbox{\sevenrm -finite}} $$
and define \hbox{$\f : A_1^{\go g} \arrow \( A_1,A_j \)$} by 
\hbox{$\f(\alpha)(\xi) = \alpha \wedge \xi$} and 
$\mu : A_1^{\go g} \t A_1^{\go g} \arrow A_3^{\go g}$ by 
$\mu(\alpha,\beta) := \alpha \wedge \beta$. {\bf QED} \bs 
{\bf Proof of the Multiplicative Eilenberg-MacLane Theorem}. 
Let $\cal C$ be the category of \hbox{$Gk$-modules}. 
Note that if $W$ is a complex in $\cal C$ then $\( k , W \)$ is the subcomplex $W^G$ 
of invariants. To any \hbox{$Gk$-module} $W$ attach the injective resolution $I(W)$ 
whose definition can be briefly recalled as follows : $I^n(W)$ is the 
\hbox{$Gk$-module} of maps from $G \times \cdots \times G$ ($n+1$ factors) to $W$, 
and the coboundary is given by
$$(df)(g_0, \dots , g_n) = \sum \ (-1)^i \ f(g_0, \dots , \widehat{g}_i, \dots, g_n).$$ 
Also remember that the (combinatorial) cup-product 
$I^p(k) \times I^q(W) \arrow I^{p+q}(W)$ is given by 
$$(\alpha \cup \beta)(g_0, \dots , g_{p+q}) = 
\alpha(g_0, \dots , g_p) \ \beta (g_p, \dots , g_{p+q}). $$
Put \hbox{$V_1 = V_2 := k$}, $V_3 := V$, $A_i := I(V_i)$ 
and define \hbox{$\f : A_1^G \arrow \( A_1,A_j \)$} by 
$\f(\alpha)(\xi) = \alpha \cup \xi$ and 
$\mu : A_1^G \t A_1^G \arrow A_3^G$ by 
$\mu(\alpha,\beta) := \alpha \cup \beta$. {\bf QED} \bs
{\bf Proof of the Multiplicative \v Cech Theorem}.
For any sheaf $T$ of $C$-modules let
$$\Psi_{T}: \check{H}^\.({\cal U},T) \to H^\.(X,T)$$
be the canonical morphism of graded $k$-modules. Consider the diagram
$$\diagram{
\check{H}^\.(X,S) \t \check{H}^\.(X,C) & 
\hfl{ \ \Psi_S \t \Phi_C \ \ }{} & 
H^\.(X,S) \t H^\.(X,C) \cr
\vfl{c}{} &                    & \vfl{}{g} \cr
\check{H}^\.(X,S)  & \hfl{}{\Phi_S} & H^\.(X,S), \cr}$$
where $c$ and $g$ denote respectively the \v Cech and the Godement cup-product as defined in section II.6.6 of [God]. Since this diagram commutes by observation (c) on page 257 of [God], it suffices to check that 
$$\diagram{
H^\.(X,S) \t H^\.(X,C) & \hfl{=}{} & H^\.(X,S) \t H^\.(X,C) \cr
\vfl{g}{} &                    & \vfl{}{t} \cr
H^\.(X,S)  & \hfl{}{=} & H^\.(X,S), \cr}$$
where, remember, $t$ is the ``true'' cup-product, commutes. Let 
$\e_C : C \to A$ and $\e_S : C \to B$ be Godement's canonical resolutions and let $\mu_A : A \t_C A \to A$, $\mu_B : B \t_C A \to B$ be the maps defined at the bottom of page 256 of [God]. To apply Propostion~1 we need maps
$$\f : \(C,A\) \to \(A,A \),$$
$$\psi : \(C,B\) \to \(A,B \),$$
$$\mu : \(C,B\) \ \t_k \ \(C,A\) \to \(C,A\).$$
Let's define them by setting
$$\big(\f(\alpha)\big)(a) := \mu_A \big(\alpha(x) \t a\big),$$
$$\big(\psi(\alpha)\big)(a) := \mu_B \big(\alpha(x) \t a \big)$$
for all $x \in X$ and all $a \in A(x)$ [the stalk over $x$] and letting $\mu$ be the composition of 
$\(C,\mu_B\)$ with the canonical map 
$$\(C,B\) \ \t_k \ \(C,A\) \to \(C,B \t_C A\).$$
Section II.6.6 of [God] implies then the assumptions of Proposition~1 are fulfilled. {\bf QED} \bs
%
%
\cl{\uq Proof \ of \ Proposition \ 1} \bs
I first claim that the diagram 
$$\diagram{
[V_i,A_j] & \hfl{H(\f)}{} & [A_i,A_j]  \cr 
\vfl{can_{{}_{V_i A_j}}}{} &           & \vfl{}{can_{{}_{A_i A_j}}} \cr 
\Ext(V_i,A_j) & \hfle{}{\Ext(\e_i,A_j)}  & \Ext(A_i,A_j), \cr} $$ 
commutes. Indeed, the restriction $[\e_i,A_j] : [A_i,A_j] \arrow [V_i,A_j]$ satisfying 
\hfill

$[\e_i,A_j] \ H(\f) = \hbox{Id}_{[V_i,A_j]}$ by assumption (a) and 
$$can_{{}_{V_i A_j}} \ [\e_i,A_j] = \Ext(\e_i,A_j) \ can_{{}_{A_i A_j}} \, ,$$
we have 
$$can_{{}_{V_i A_j}} = can_{{}_{V_i A_j}} \ [\e_i,A_j] \ H(\f) 
= \Ext(\e_i,A_j) \ can_{{}_{A_i A_j}} \ H(\f).$$ 
The top square of
$$\diagram{
[V_2,A_3] \ \t \ [V_1,A_2] \ & \hfl{\mu_*}{} & [V_1,A_3]  \cr
\vfl{H(\f_{{}_{23}}) \t H(\f_{{}_{12}})}{} & & \vfl{}{H(\f_{{}_{13}})} \cr 
[A_2,A_3] \ \t \ [A_1,A_2]                 & \hfl{c}{}  & [ A_1,A_3]  \cr
\vfl{can_{23} \t can_{12}}{}           &             & \vfl{}{can_{13}} \cr 
\Ext(V_2,V_3) \ \t \ \Ext(V_1,V_2)                 & \hfl{}{c}  & \Ext(V_1,V_3)
\cr} $$ 
commutes by assumption (b)~; the bottom square commutes because it is of the form (1)~; 
the vertical compositions are respectively $can_{23} \t can_{12}$ and $can_{13}$ 
by the claim. {\bf QED} \bs
\cl{\uq Proof \ of \ Strong \ Proposition \ 2} \bs
For any complex $C$ let $C[n]$ be the 
complex $C^{n + \.}$~; for any complex morphism \hbox{$f : B \arrow C[n]$} denote by 
$[f] \in [B,C]^n$ its homotopy class and by 
$\wt f$ the corresponding element of $\Ext^n(B,C)$~; any object of $\cal C$ shall be 
viewed as a complex in degree zero. --- Let $V$ be an object and $A$ 
a complex satisfying $A^n = 0$ for $n < 0$ and $H^n(A) = 0$ for $n > 0$~; let 
$Z \subset A$ be the subcomplex of cocycles (in other words $A$ is a right resolution 
of $Z^0$)~; set 
$$F^n := \Ext^n(V,-).$$
By left exactness of $F^0$ the canonical 
morphism $H^0(F^0 A) \arrow F^0 Z^0$ is an isomorphism, proving Strong Proposition 2 in 
degree 0. For $p > 0$ the short exact 
sequence \hbox{$Z^{p-1} \mono A^{p-1} \epi A^p$} gives birth to the 
long exact sequence
$$\diagram{
 &  &  &  & 0 & \hfl{}{} \cr \cr
F^0 Z^{p-1} & \hfl{}{}  & F^0 A^{p-1} & \hfl{}{} & F^0 Z^p & \hfl{\delta_{p,0}}{} \cr \cr
F^1 Z^{p-1} & \hfl{}{} & F^1 A^{p-1} & \hfl{}{} & F^1 Z^p & \hfl{\delta_{p,1}}{} \cr \cr
F^2 Z^{p-1} & \hfl{}{} & F^2 A^{p-1} & \hfl{}{} & \cdots      & \cdots \cr \cr
\vdots  & \vdots & \vdots  & \vdots & \vdots      & \vdots   \cr}$$
(see [V,III.1.2.5,p.162]). 
Let $\overline \delta_{p,0} : H^p(F^0 A) \arrow F^1 Z^{p-1}$ be the monomorphism induced by 
$\delta_{p,0}$ and for $r \geq 0$ introduce the inclusion
$$i_r : Z^r \incl A[r+1].$$
Fix a positive integer $n$. We want to interpret in terms of connecting morphisms the 
canonical morphism $H^n(F^0 A) \arrow \Ext^n(V,Z^0)$, which I prefer to think as the 
morphism 
$$can : H^n(F^0 A) \arrow \Ext^n(V,A)$$
satisfying 
$$can([x]) = \wt i_n \ \wt x \eqno{(2)}$$
for all morphisms $x : V \arrow Z^n$. Define 
$$\psi : H^n(F^0 A) \arrow \Ext^n(V,A)$$
by 
$$\psi := \Ext(V,i_0) \ \delta_{1,n-1} \ \delta_{2,n-2} \ \cdots \ 
\delta_{n-1,1} \ \overline \delta_{n,0} \ .$$
Strong Proposition 2 follows easily from \bs
{\bf Lemma}. We have $can = (-1)^{n(n+1)/2} \ \psi$. \bs
{\bf Proof}. For $p > 0$ let $\iota$ be the inclusion of $Z^{p-1}$ into $A^{p-1}$ 
and $C_p$ the complex 
$$\diagram{C_p^{-1} = Z^{p-1} & \hfl{- \iota}{} & A^{p-1} = C_p^0 \ , \cr}$$
and consider the morphisms 
$$\diagram{
Z^p & \hfle{q_p}{} & C_p & \hfl{f_p}{} & Z^{p-1}[1]  \cr}$$
defined by
$$\diagram{
Z^p & \hfle{d}{} & A^{p-1} & & \cr
           &  & \vfle{- \iota}{} & \cr
 &   & Z^{p-1} & \hfl{=}{} & Z^{p-1}[1] \cr \cr}$$
Note the following : \hbox{$\wt f_p \in \Ext^1(C_p,Z^{p-1})$~;}  
\hbox{$\wt q_p \in \Ext^0(C_p,Z^p)$~;} $\wt q_p$ is a quasi-isomorphism~;  
\hbox{$\wt f_p \ \wt q_p^{\, -1} \in \Ext^1(Z^p,Z^{p-1})$.} Recall 
that $\delta_{p,r} : F^r Z^p \arrow F^{r+1} Z^{p-1}$ coincides with the left 
multiplication by $\wt f_p \ \wt q_p^{ \, -1}$ (see [Iv,XI.3])~; in particular 
$$\psi([x]) := \wt{i_0} \ \wt f_1 \ \wt q_1^{ \, -1} \, \wt f_2 \ \wt q_2^{ \, -1} \ 
\cdots  \ \wt f_{n-1} \ \wt q_{n-1}^{ \, -1} \, \wt f_n \ \wt q_n^{ \, -1} \,
\wt x \eqno{(3)}$$
for $x : V \arrow Z^n$. Confronting (2) and (3) we see that the lemma reduces to the 
equality 
$$\wt{i_0} \ \wt f_1 \ \wt q_1^{ \, -1} \ \wt f_2 \ \wt q_2^{ \, -1} \ 
\cdots  \ \wt f_{n-1} \ \wt q_{n-1}^{ \, -1} \  \wt f_n \ \wt q_n^{ \, -1}
 = (-1)^n \ \wt i_n \, .$$
It suffices thus to check that for $p > 0$ we have  
$\wt i_{p-1} \ \wt f_p \ \wt q_p^{ \, -1} = (-1)^p \ \wt i_p \, ,$
that is
$\wt i_{p-1} \ \wt f_p = - \ \wt i_p \ \wt q_p \, .$
The morphisms $i_{p-1} \, f_p$ and $i_p \, q_p$ from $C_p$ to $A[p]$ being respectively 
given by the diagrams
$$\diagram{ & & \vdots \cr
 & & \vfle{}{(-1)^p d} \cr
A^{p-1} &  & A^p \cr
\vfle{- \iota}{} &  & \vfle{}{(-1)^p d} \cr
Z^{p-1} & \hfl{\iota}{} & A^{p-1} \cr
 & & \vfle{}{(-1)^p d} \cr
 & & \vdots \cr}$$
and
$$\diagram{ & & \vdots \cr
 & & \vfle{}{(-1)^p d} \cr
A^{p-1} & \hfl{d}{} & A^p \cr
\vfle{- \iota}{} &  & \vfle{}{(-1)^p d} \cr
Z^{p-1} & & A^{p-1} \cr
 & & \vfle{}{(-1)^p d} \cr
 & & \ \vdots \ , \cr}$$
the identity of $A^{p-1}$ furnishes a homotopy from $i_{p-1} f_p$ to 
$(-1)^p \ i_p \, q_p$. {\bf QED} \bs 
\cl{* \ * \ *} \bs
{\parindent=11mm 
\item{[BT]} Bott R. \& Tu Loring W., {\bf Differential forms in algebraic topology}, 
Graduate Texts in Mathematics, 82, Springer-Verlag. XIV, 1982. \bs
\item{[BW]} Borel A.\& Wallach N., {\bf Continuous cohomology, discrete subgroups, and 
representations of reductive groups}, Annals of Mathematics Studies, No.~94, 
Princeton University Press 1980. \bs
\item{[CE]} Cartan H. \& Eilenberg S., {\bf Homological algebra}, Princeton 
University 
Press 1956. \bs
\item{[G]} Grothendieck A., Sur quelques points d'alg\`ebre homologique, 
{\it T\^ohoku Math. J.} {\bf 9} (1957) 119-221. \bs
\item{[God]}Godement R., {\bf Topologie algebrique et theorie des faisceaux}, Hermann, 
Paris 1958.
\item{[Iv]} Iversen B., {\bf Cohomology of sheaves}, Springer 1986. \bs
\item{[V]} Verdier, J.-L., {\bf Des Cat\'egories d\'eriv\'ees des cat\'egories 
ab\'eliennes}, Ast\'erisque 239, SMF, Paris 1996. \bs 
\item{[W]} Warner, F., {\bf Foundations of differentiable manifolds and Lie groups}, 
Springer-Verlag 1987. 
\par} \bs
\hfill Last update~: Feb. 27, 2000 \bs
This text and others are available at http://www.iecn.u-nancy.fr/$\sim$gaillard
\bye